\documentclass[conference]{IEEEtran}
\usepackage[noadjust]{cite}
\usepackage{graphicx,color,psfrag,subeqnarray}
\usepackage{pgfplots}
\usepgfplotslibrary{patchplots,colormaps}
\usepackage{amsmath}
\usepackage{amssymb}
\usepackage{amsthm,amssymb}
\usepackage{color}
\usepackage{listings}
\usepackage{multirow}
\usepackage{bm,array}
\usepackage{listings}
\usepackage[font=footnotesize]{caption}
\usepackage{subcaption}
\usepackage{float} 
\usepackage{bigstrut}
\usepackage{epstopdf}

\usepackage{algorithm}
\usepackage[noend]{algpseudocode}

\makeatletter
\def\BState{\State\hskip-\ALG@thistlm}
\makeatother

\DeclareCaptionSubType[alph]{figure}
\newcommand\textlcsc[1]{\textsc{\MakeLowercase{#1}}}

\begin{document}

\title{Decarbonized Demand Response for Residential Plug-in Electric Vehicles in Smart Grids}

\author{\IEEEauthorblockN{Farshad Rassaei, Wee-Seng Soh and Kee-Chaing Chua  
		 \\}
\IEEEauthorblockA{Department of Electrical and Computer Engineering\\
National University of Singapore, Singapore\\
Email: farshad@u.nus.edu, \{weeseng,chuakc\}@nus.edu.sg}
}
\maketitle

\begin{abstract}
Recently, in Paris, the world has reached an agreement whereby many countries commit to bolster their efforts about reducing adverse climate changes. Hence, we can expect that \textit{decarbonization} will even attract more attention in different energy sectors in near future. In particular, both generation side and consumption side are required to be run more congruently and environmentally friendly. Thus, employing the renewables at the generation side along with our proposed \textit{decarbonized demand response} (DDR) at the consumption side could significantly reduce deleterious impacts on the climate.
Such ambition, at the consumption side, necessitates \textit{symbiosis} and \textit{synergy} between the customers and the retailer, and among customers, respectively. In other words, there should be some incentive-based collaboration between customers and the retailer as well as coordination among customers to make the objective be achieved successfully. In this paper, we present such matching demand response (DR) algorithm for residential users owning vehicle-to-grid (V2G) enabled plug-in electric vehicles (PEVs) who obtain electricity from a common retailer. The retailer itself is connected to the wholesale electricity market to purchase and sell electricity. Furthermore, we explain the details of the existing symbiosis and synergy in our system. Our simulation results illustrate that substantial cost savings can be achieved along with pollution reduction by our proposed technique.              
\end{abstract}


\IEEEpeerreviewmaketitle

\section{Introduction}
Climate change has become one of the major concerns worldwide. Recently, in December 2015, many countries have agreed to further enhance their efforts to confront adverse climate changes which are mainly because of tremendous green house gases (GHGs) emissions, e.g., $\text{CO}_2$ and $\text{CH}_4$ \cite{burleson2016paris}. 

One of the significant reasons for GHG emissions is the transportation sector. Thus, decarbonization this sector has attracted much research, e.g, \cite{zhang2016times,WENE:WENE181,liu2015decarbonizing,abdallah2015,mittal2015low}. Meanwhile, plug-in electric vehicles (PEVs) are good alternatives for traditional cars to diminish carbon emissions \textit{provided} their electricity consumption is managed properly.    

However, in the literature, PEVs' charging and/or discharging management and scheduling are mainly investigated for cost savings purposes, e.g., \cite{kim_bidirectional_2013,mohsenian-rad_autonomous_2010,sgtrfar,jointfar2}. In these papers, the emphasis is mostly on increasing the users' utility, welfare, the billing strategies, etc. 

On the other hand, in near future, we are going to face a new paradigm in power system, e.g., new ways of electricity generation, market liberalization, storage capability, two-way electricity delivery, demand side management (DSM), demand response (DR) and environmentally concious transportation \cite{marnay2013challenges} and \cite{kutz2008environmentally}. Let us add the salience of decarbonization to the above list. Hence, a practical technique is incumbent to consider this new paradigm in order to be competent enough to be employed in a real-world smart grid. 

\begin{figure}[t]
	\centering
	\includegraphics[width=\columnwidth,height=3in]{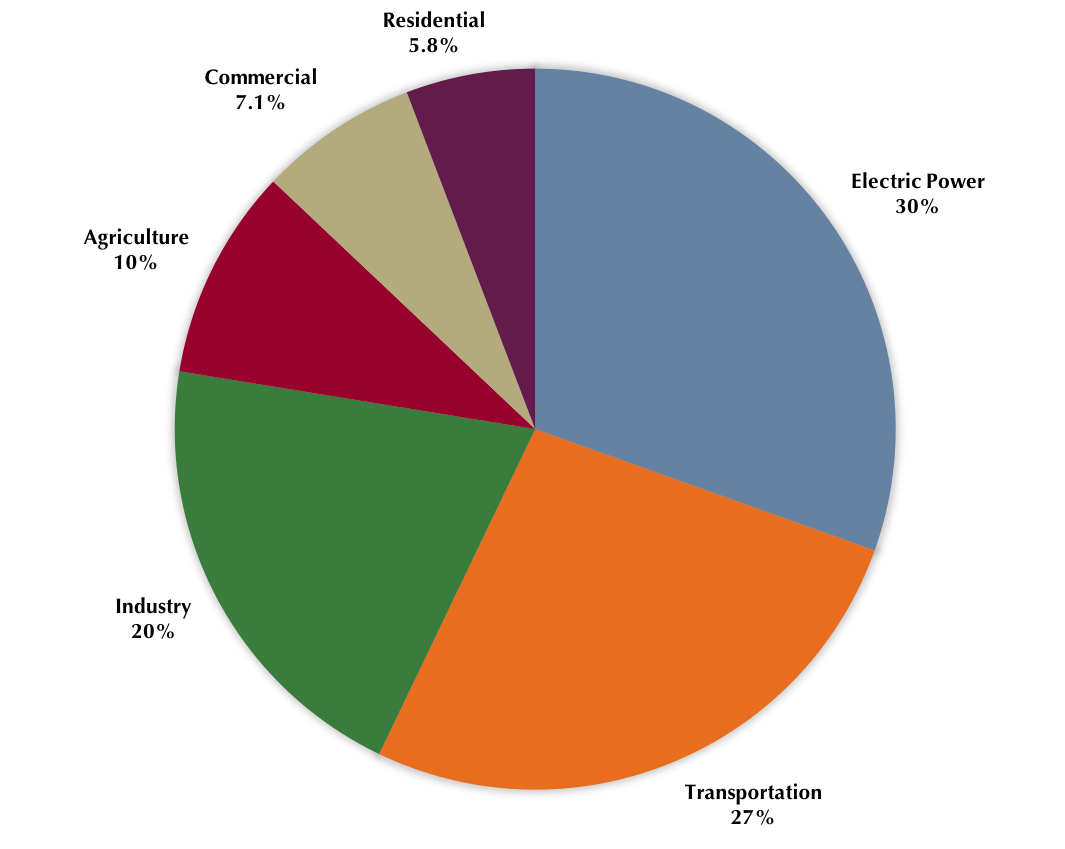}
	\caption[]{The United States of America GHGs emissions by sector in 2014 \cite{ce14}.}	 
	\label{ce}\vspace{-1.5 em}
\end{figure}

The share of the electricity generation and transportation sectors in GHGs emissions depend on many factors, e.g., type and age of generators, different regions, traffic congestion management, etc. It varies from one region to another one in the world. For instance, Fig. \ref{ce} illustrates the GHGs emissions by sector in 2014 in the U.S. and we can observe that the quota from electricity generation and transportation sectors stand for more than half of the overall emissions.  

\begin{figure}[t]
	\centering
	\includegraphics[width=\columnwidth,height=2.2in]{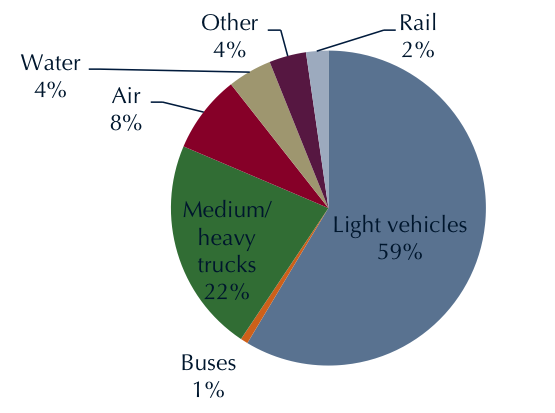}
	\caption[]{The share of GHGs emissions in the transportation sector by mode in the U.S. \cite{ce14}.}	 
	\label{dift}
\end{figure}

However, we believe that minimizing and even nullifying the share of GHGs emissions from the transportation and the power generation sectors will be possible in the long run by further utilizing renewables, electromobility and proper DDR techniques. Besides, the emissions from power generators could be reduced in smart grids wherein there is high penetration of renewables and distributed generation (DG). Widespread penetration of PEVs could inherently reduce the GHGs emissions from the transportation sector. Fig. \ref{dift} shows that the share from light duty vehicles (LDV) accounts for 59\%. Nonetheless, PEVs' electricity demand adds a huge burden on the power generation side. 

We should note that striving to make decarbonized energy supply alone is not adequate \cite{williams2012technology} nor is electrification of transportation sector. In order to triumph in GHGs reduction, congruous DR techniques are also needed. 
In other words, taking into account the level of GHGs emissions for the generation and consumption sides in DR can further diminish the emissions even in an already electrificated sector. 

\begin{figure} [t]
	\centering
	\includegraphics[width=\columnwidth,height=2.4in]{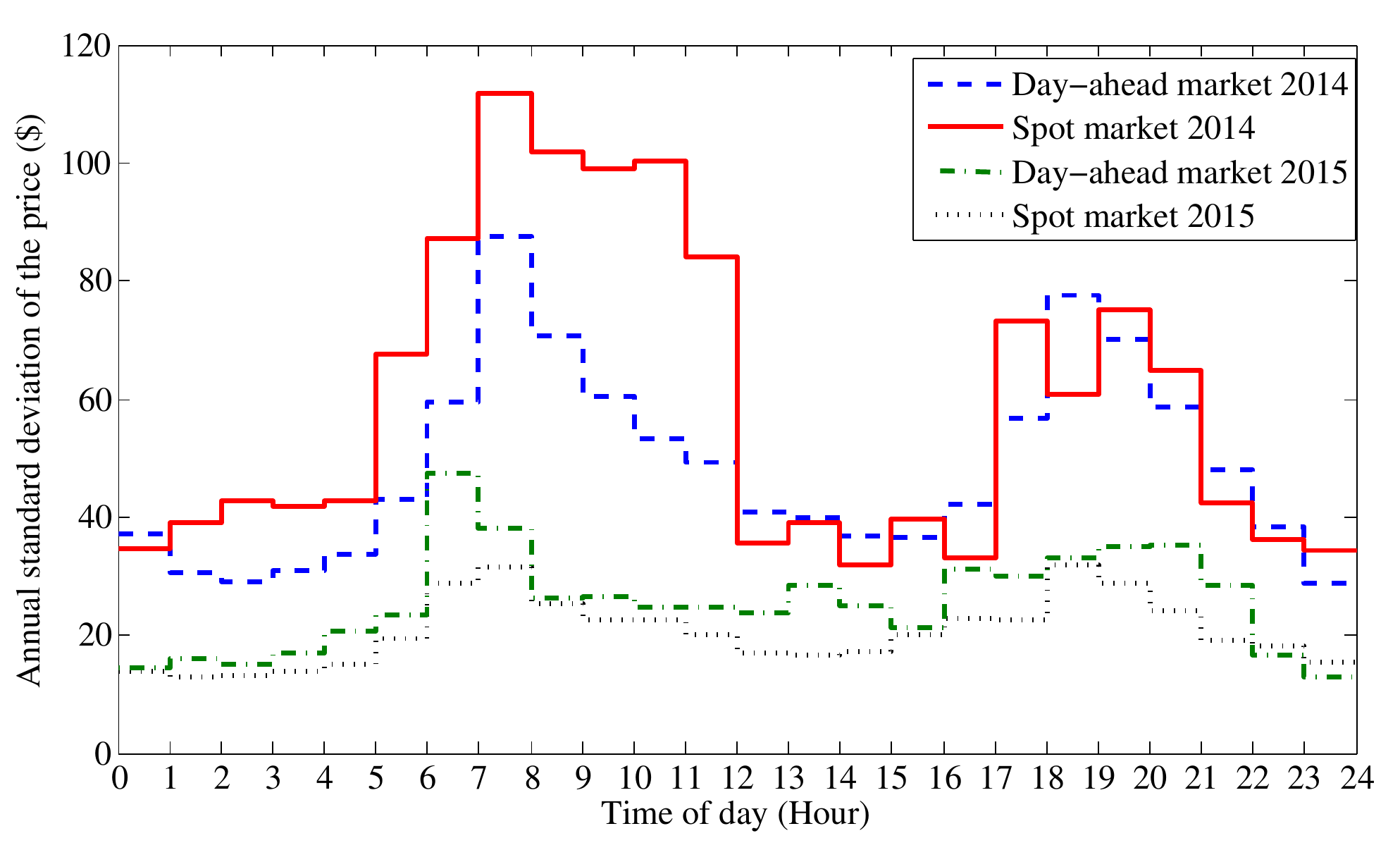}
	\caption[]{Annual standard deviation of the electricity price in 2014 and 2015 for each our of a day in PJM for DA and RT (spot) markets.}	 
	\label{std45}
\end{figure}

%In a deregulated electricity system, electricity retailers submit demand bids to wholesale electricity markets. For instance, for a day-ahead (DA) electricity market, these demand bids often have both power demand's amounts and price components meaning that the retailer purchases the specified power amount only if the market clearing price (MCP) is not more than the price component \cite{mohsenianoptimal}. This bidding can be carried out in a few predefined rounds allowing the retailers to update their bids at each round before final clearance. This type of bidding is referred to as limit order bidding. Therefore, in this case, it is admissible to presume that a retailer is willing to shape its aggregated demand profile and match it to the electricity profile resulted from successful bids for different hours of the following day so that it can minimize its demand from the RT market to balance the load, and accordingly reduce the overall electricity procurement cost for each following day.  

In this paper, we consider  Pennsylvania–-Jersey-–Maryland (PJM) day-ahead (DA) and real-time (RT) electricity market. We use its pricing data for the years 2014 and 2015 \cite{PJM}, the average prices of electricity per MWh which has been sold over those two years are very close for both DA and RT markets in each year, i.e., \$48.95 and \$48.21, in 2014, and \$33.94 and \$33.34 in 2015, respectively. The reason for cheaper average price in 2015 compared to its antecedent year could be the unprecedented falling down of the oil price. 

Fig. \ref{std45} shows the annual standard deviation of electricity price in 2014 and 2015 for each hour of a day in PJM for both DA and RT (spot) markets. Although the prices are much cheaper in 2015, we observe that hourly pricing data for PJM's DA and RT markets can be significantly distinct and unpredictable. Another point is that we see RT market prices has more fluctuations than the DA market, as we could expect it. Therefore, the high uncertainty, particularly in the RT market, can remarkably affect the overall electricity procurement cost for a retailer especially in the long term. This fact is much more expected when the power system is relying on a large number of intermittent energy resources with more uncertainty.

The role of intelligence along with significant architectures and concepts in future power systems are reviewed in \cite{strasser2015review}. A good overview of DR and their different classifications in a deregulated electricity market is discussed in \cite{albadi2008summary}. In \cite{ISFAR}, we present a statistical modelling and a closed-form expression for PEVs' uncoordinated charging power demand. Our paper \cite{sgtrfar} proposes a decentralized demand shaping algorithm for an a priori known demand profile for the next day or for flattening the aggregated daily demand profile. In \cite{jointfar2} we consider both DA and RT markets of PJM in our proposed DR algorithm.

In this paper, for the transportation sector, by adding the significance of reducing GHGs emissions, we discuss \textit{decarbonized} DR (DDR) techniques for residential users owning vehicle to grid (V2G) enabled PEVs by which we strive to decrease the emissions from the electric power sector, see Fig. \ref{ce}. Hence, we contemplate lessening both carbon emissions and electricity procurement costs.

\begin{figure}[t]
	\centering
	\includegraphics[width=\linewidth,height=2.7in]{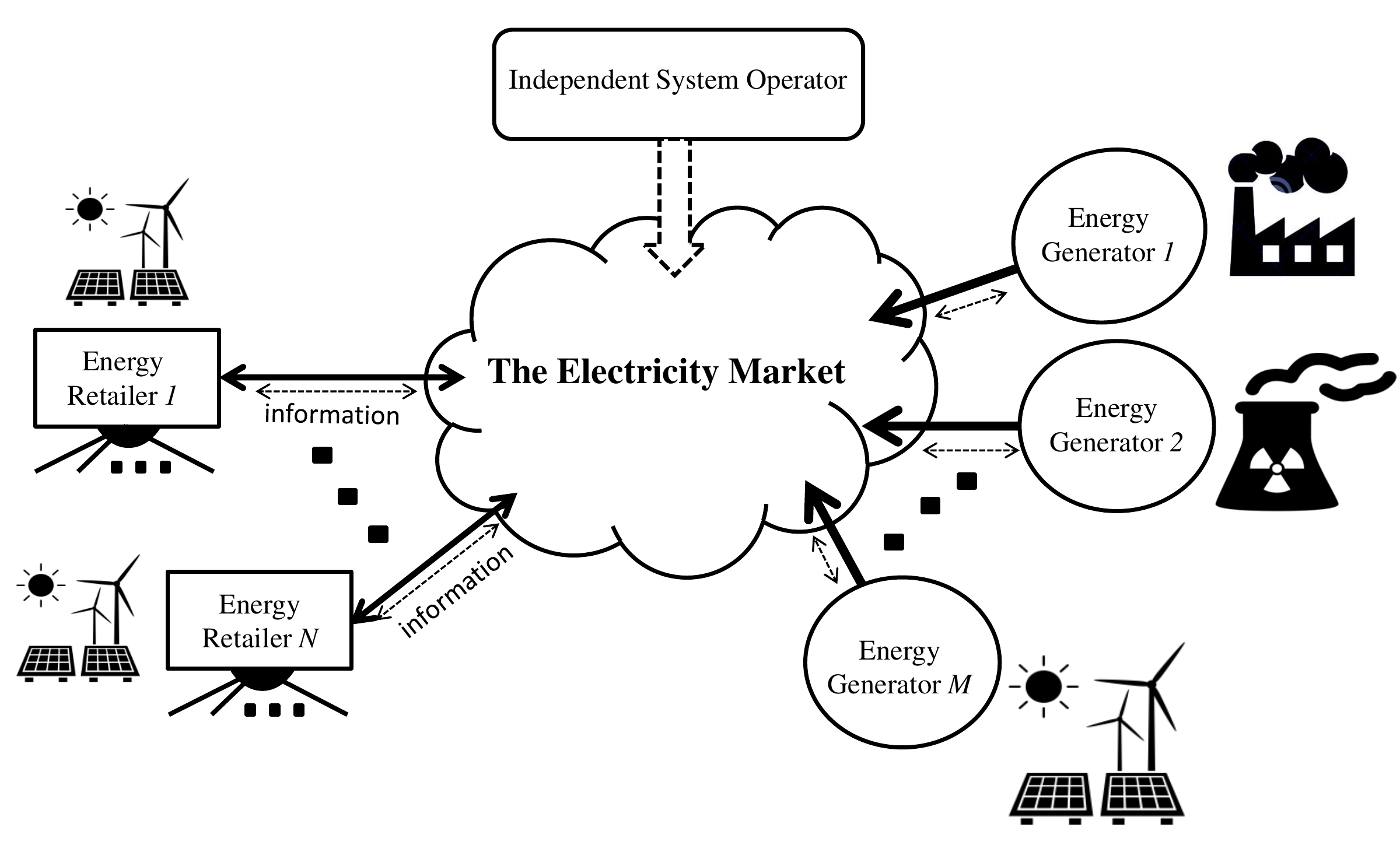} 	
	\caption{An archetype of a smart electricity system encompassing several retailers, utility companies and multiple users sharing one electricity retailer or an aggregator.} 
	\label{mod} 
\end{figure}

\section{System Model} \label{SM}
In this section, we provide the underlying model and assumptions of the power system in this paper which entails the energy markets, the electricity retailers or the aggregators, and the residential users. Similar models for future smart power systems are advocated in \cite{mohsenian-rad_autonomous_2010} and \cite{liu2014peak}. We discuss this model in the sequel.

Fig. \ref{mod} represents our model of a smart electricity system where multiple users share one electricity retailer or an aggregator. The users' overall load can be differentiated into two distinct types; typical household load which normally needs \textit{on-demand} electricity supply, e.g. air conditioning, heating, lighting, audio visual devices, cooking and refrigerator, and PEV as a \textit{flexible} or programmable load. Here, the dotted lines show the underlying information system and the solid lines represent the power transmission and distribution infrastructure. 
 
We assume that an electricity retailer (which may own its generation capacity) bids to the energy market, e.g., on a DA basis. Then, based on its energy needs and the market situation, it buys electricity from the market at market clearing prices (MCPs). Then, we assume that the retailer is willing to handle its customers' PEVs' electricity assignments (charging and discharging) such that the shape of the resulting aggregated power demand profile matches the electricity profile resulted from the successful bids in the DA market. 

This enables the retailer to minimize its demand from the RT market -which has more price volatility according to Fig. \ref{std45}- for balancing the load in the following day. Accordingly, it can reduce the overall electricity procurement cost. This cost reduction makes the energy retailer afford to offer more attractive deals to the customers in the form of pricing, rewarding, promotions, etc \cite{sgtrfar}.

On the other hand, we assume that there are some incentives or limits from a regulator or the government which make the retailers interested in or have to reduce the GHGs emissions. We note that the incentives and limits can be translated to payoffs and fines, in terms of fulfilment and violation, respectively. As we indicated in \cite{jointfar2}, retailers adjust their electricity deals (purchase and sell) in response to market prices. Nevertheless, the regulator can put some limits on power consumption, e.g., the mean of daily electricity consumption. This could reduce or cancel the need of turning on the traditional generators which are accounted for 30\% of the GHGs emissions Fig. \ref{ce}. 

In practice, the shaped aggregated profile does not exactly suite the retailer's purchased DA energy profile. Thus, the retailer needs to reciprocate the load imbalances in the following day by referring to the RT market and buy electricity at RT prices. Therefore, we assume that a retailer should consider three directions when designing its DR technique: DA market, RT market and GHGs emissions. 

%On the other hand, we should notice that residential DR is desired to be carried out such that users' privacy is not violated. In such case, we can expect that DR attracts more participation from the users, especially when the main job is implemented in each user's house in a decentralized fashion according to this model.   

\section{Analysis} \label{Sta}
In this section, we provide the electricity procurement cost for a retailer bidding to the electricity market and then present our proposed DDR algorithm.

%\hspace{-1.5em}\begin{tikzpicture}
%\begin{axis}[%
%width=\columnwidth,height=8.5cm,
%xlabel={$P_1$},
%ylabel={$P_2$},
%zlabel={$P_3$},
%legend style={at={(-0.2,0.14)},anchor=north west,draw=black,fill=white,legend cell align=left},
%label style={font=\scriptsize}, ticklabel style={font=\scriptsize}
%]
%
%\addplot3[patch,patch type=triangle quadr,opacity = 0.2,
%shader=faceted interp]
%coordinates {
%	(349.9671,  349.9671,  195.8676)
%	(195.8676,  349.9671,  349.9671)
%	(349.9671,  195.8676,  349.9671)
%	(226.6197,  330.3199,  226.6197)
%	(226.6197,  226.6197,  330.3199)
%	(330.3199,  226.6197,  226.6197)};
%\end{axis}
%\end{tikzpicture}
%\begin{figure}[h]
%	\caption{An example of Pareto-optimal surface for the multi-objective programming \textbf{P}.} 
%	\label{pareto} 
%\end{figure} 

The overall electricity procurement cost for the next day can be formulated as follows: 
\begin{align}
\nonumber Cost=<\textit{\textbf{l}}^{DA},\textit{\textbf{p}}^{DA}>+<\textit{\textbf{l}}^i,\textit{\textbf{p}}^{RT}>\\
=
\sum_{t=1}^{24}p^{DA}_{t}l^{DA}_{t}+\sum_{t=1}^{24}p^{RT}_{t}l^i_{t},
%=<\textit{\textbf{l}}^d,\textit{\textbf{p}}^d>+<\textit{\textbf{l}}^i,\textit{\textbf{p}}^r>,
\end{align}

\noindent where $Cost$ is the overall electricity procurement cost over the scheduling horizon. The 24-elements vectors $\textit{\textbf{l}}^{DA}$ and $\textit{\textbf{p}}^{DA}$ represent the power demand and electricity price cleared in the DA market for the following day. The units of $l_i^{DA}$ and $p_i^{DA}$ are MWh and \$/MWh, respectively. Similarly, $\textit{\textbf{l}}^i$ and $\textit{\textbf{p}}^{RT}$ are load imbalance and the electricity price vectors in the RT market for the next day. The values of the elements of these two vectors will be known to the retailer only at each time slot of the next day.

%\vspace{-2em}

Then, the following sequential optimization programming technique is used. The users individually contribute in this program (see Algorithm 1): 

\begin{align}
\nonumber \underset{\textit{\textbf{l}}'_{\text{PEV},n}}
 {\text{minimize}} \quad
(\lambda <\textit{\textbf{l}}'_{\text{PEV,n}},\textit{\textbf{l}}'_{\text{A},n} + \textit{\textbf{l}}'_{-n}-\textit{\textbf{l}}^d>\\ 
+(1-\lambda) (l'^{t_0}_{-n}+l^{t_0}_{\text{A},n}+l'^{t_0}_{\text{PEV},n})),
\label{P}
\end{align}\vspace{-1.2em}
\begin{align}
&[l^{'1}_{\text{PEV},n},\cdots,l^{'t_0-1}_{\text{PEV},n}]=[l^{1}_{\text{PEV},n},\cdots,l^{t_0-1}_{\text{PEV},n}]   ,\\
&\sum_{t=t_0}^{\beta_{n}} l'^{t}_{\text{PEV},n}=E_{\text{PEV},n}-\sum_{t=\alpha_{n}}^{t_0-1} l^{t}_{\text{PEV},n},  \\ 
&\sum_{n=1}^{N}(l^{t}_{\text{PEV},n}+l^{t}_{\text{A},n})=1.5 \times \frac{\sum_{t=1}^{24}\sum_{n=1}^{N}(l^{t}_{\text{PEV},n}+l^{t}_{\text{A},n})}{N},  \quad \forall t,\\
%& l^{t}_{\text{PEV},n}=l^{t}_{\text{PEV},n}, \quad t\\
&|l'^{t}_{\text{PEV},n}| \leq p_{max},\quad \forall t\in\mathbb{T}^P_{\text{PEV},n}, \\
%& SOC^t_{\text{PEV},n}\geq 0.2\times C_{\text{PEV},n}, \quad \forall t\in\mathbb{T}^P_{\text{PEV},n}, \\ 
& l'^{t}_{\text{PEV},n}=0, \quad \forall t\notin \mathbb{T}^P_{\text{PEV},n},\\ 
& SOC^{t=t_0-1}_{\text{PEV},n}+\sum_{k=t_0}^{t}l'^k_{\text{PEV},n}\geq 0.2 \times C_{\text{PEV},n}, 
\forall t\in\mathbb{T}^P_{\text{PEV},n},
\end{align}
\noindent Here, $\textit{\textbf{l}}_{\text{PEV},n}$ and $\textit{\textbf{l}}_{\text{A},n}$ show the energy assignment vectors for user n's PEV and the aggregated load from its household appliances, respectively. Besides, $\textit{\textbf{l}}'_{\text{PEV},n}$ and $\textit{\textbf{l}}'_{\text{A},n}$ vectors show the same things whenever load altering is needed in RT market, see Algorithm 1. Furthermore, $\textit{\textbf{l}}^d$ is the purchased load profile from DA market. $E_{\text{PEV},n}$ is the user $n$'s required energy to be delivered to its PEV which is associated with the total required charging time $T_{\text{PEV},n}$ as follows: 
\begin{align}
E_{\text{PEV},n}=\omega \times T_{\text{PEV},n},
\end{align} 
\noindent where $\omega$ is the charging power rate of the outlet to which it is plugged in. Likewise, $\alpha_{n}$ and $\beta_{n}$ represent the arrival time and departure time of the PEV. Furthermore, $|l^{t}_{\text{PEV},n}| \leq p_{max}$ limits the maximum power that can be delivered to/from the PEV, we may presume $p_{max}=\omega$, and $\mathbb{T}^P_{\text{PEV},n}$ describes the permissible charging time set or simply the set of time slots during the PEV's \textit{connection time} to the power grid. This is simply the set of time slots between $\alpha_{n}$ and $\beta_{n}$. Constraint (5) limits excessive power consumption at each time slot. This prevents the need to turn on or use traditional thermal generators.

\begin{figure}[t]
	\centering
	\includegraphics[width=\columnwidth,height=2.8in]{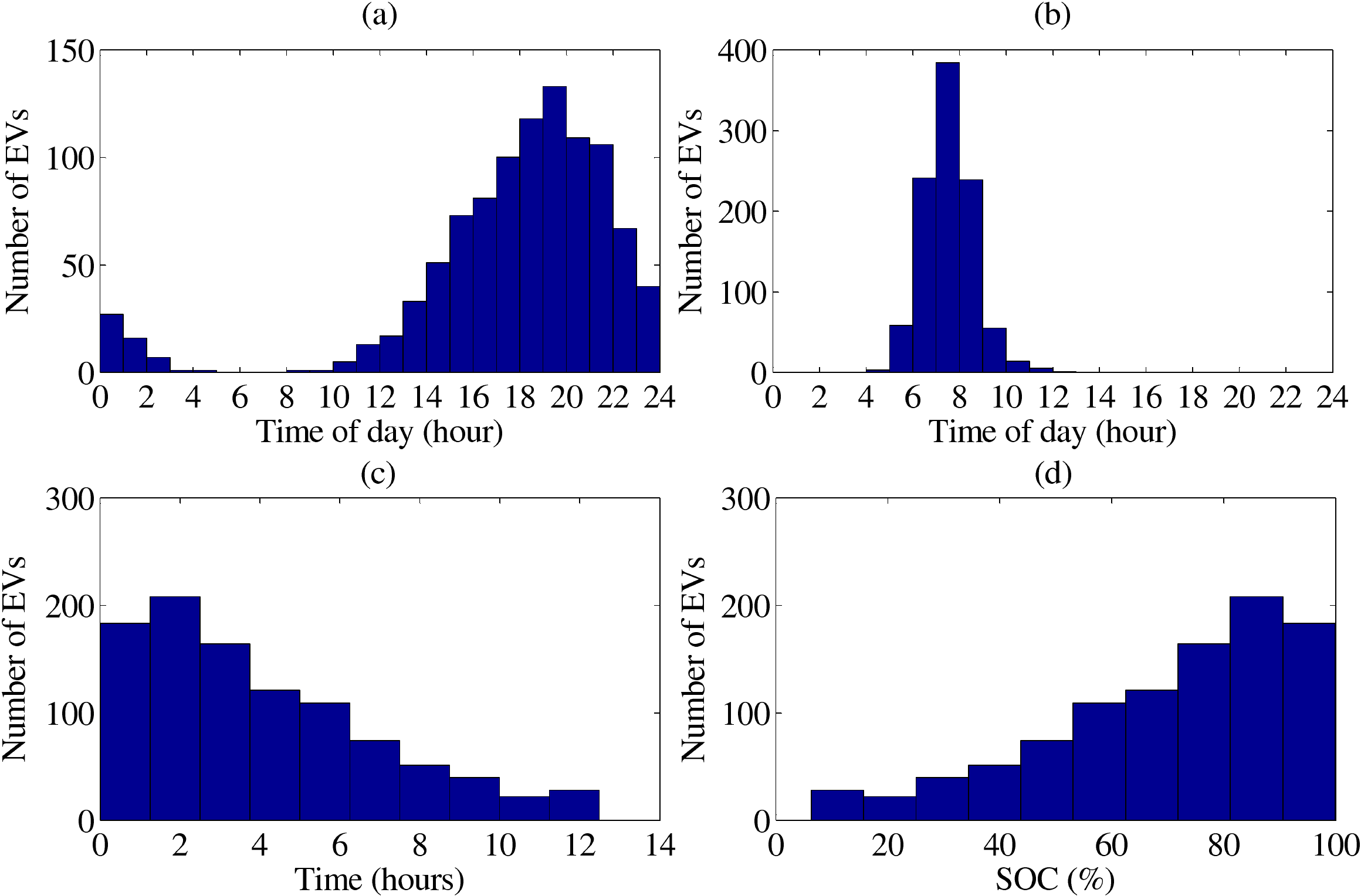}
	\caption{Distributions of (a) arrival time, (b) departure time, (c) charging time and (d) initial SOC, for 1,000 electric vehicles.}
	\label{tdist}\vspace{-2em}
\end{figure}

\begin{algorithm}[t]
	\begin{algorithmic}[1]
		\caption{Decarbonized Demand Response (DDR)}\label{DDR}
		\State Each user initializes its respective load profile over the assignment horizon based on its power demands, i.e., $\textit{\textbf{l}}_{n}$ for $n=1,\dotsc, N$.
		\State All $N$ users send their initialized load profiles to the retailer.
		\While {\textit{not reaching convergence}}
		%\State The retailer does a random ordering of the users.
		\For {$n=1$ to $N$}
		\State $\lambda$ is set $1$ in the proposed problem. 
		\State  The retailer calculates the state information $\textit{\textbf{l}}_{-n}$    \hspace*{10mm} according to (\ref{others}) for user $n$.
		\State The retailer sends $(\textit{\textbf{l}}_{-n}-\textit{\textbf{l}}^d)$ to user $n$.
		\State Each user $n$ solves the proposed problem and \hspace*{10mm} updates its load profile $\textit{\textbf{l}}_{n}$.
		\State User $n$ sends back the new demand profile to the \hspace*{10mm} retailer. 
		\State The retailer updates $\textit{\textbf{l}}_{n}$.
		\EndFor 
		\State \textbf{end for}
		%\State $i \gets i+1$
		\EndWhile 
		\State \textbf{end while}
		\For {$t=1$ to $24$}
		\State The retailer receives information from RT market, i.e., \hspace*{5mm} $\textit{p}^{RT}_t$.
		\State The retailer decides whether or not it proceeds for \hspace*{5mm} demand altering.
		\State $\lambda$ is set to a desired value in the proposed problem. 
		\While {\textit{not reaching convergence}}
		\For {$n=1$ to $N$}
		\State The retailer sends demand altering signal at \hspace*{15mm} time slot $t$ to user $n$.
		\State User $n$ solves problem (2) and updates its \hspace*{15mm} load profile $\textit{\textbf{l}}'_{n}$.
		\State User $n$ sends back the new demand profile \hspace*{15mm} to the retailer. 
		\State The retailer updates $\textit{\textbf{l}}'_{n}$.
		\State  The retailer calculates the state information \hspace*{15mm} $\textit{\textbf{l}}'_{-n}$ according to (\ref{others}) for user $n$.
		%		\State The retailer sends $\textbf{l}^{'}_{-n}$ to user $n+1$.
		\EndFor 
		\State \textbf{end for}
		\EndWhile 
		\State \textbf{end while}
		\EndFor 
		\State \textbf{end for}
	\end{algorithmic}
\end{algorithm}

Additionally, $\textit{\textbf{l}}_{-n}$ is the aggregated power profile from other $N-1$ users in the system which is described as follows: 
\begin{align}
\textit{\textbf{l}}_{-n}=\sum_{\substack{i \in {N} \\  i \neq n}}  (\textit{\textbf{l}}_{\text{PEV},i}+\textit{\textbf{l}}_{\text{A},i}). 
\label{others}
\end{align}
In (8), $C_{\text{PEV},n}$ is the total storage capacity of the user n's PEV and we assume that in case of employing V2G in the system, PEV's state of charge (SOC) should not fall below 20\% of that total capacity in order to make sure that the adverse impacts on PEV's battery lifetime due to complete depletion are avoided.   
 
We know the fact that $(\textit{\textbf{l}}^{i},\textit{\textbf{p}}^{RT})$ is unknown to the retailer a priori, at each time slot $t_0$ of the next day after getting this information, the retailer may decide to alter the previously shaped DA demand profile. It may want to minimize its RT electricity purchase to balance the load if the RT prices rise unexpectedly and even sell back some of its pre-purchased electricity from DA market to the RT market by using the PEVs' available demand elasticity. RT prices may fluctuate significantly due to the state of the RT market or contingencies.

We should notice that in the proposed programming method $\lambda$ is one for shaping the aggregated demand profile, $\textit{\textbf{l}}_N$, for the DA market. Then, as described in Algorithm 1, $\lambda$ could take a desired value to alter the aggregated demand profile, $\textit{\textbf{l}}'_N$, in response to price fluctuations in the RT market.  

As the chances for the price to remain that high during all the next remaining hours of the day is low \cite{poterba1988mean}, reshaping the load profile by lowering the electricity consumption at that time slot and purchasing electricity at the further time slots can yield to a lower electricity procurement total cost in practice. This is also true for purchasing electricity at those time slots when price, unexpectedly, falls down significantly. The retailer may buy extra electricity at those specific time slots (based on the overall storage capacity coming from connected PEVs).

We should note that the retailer is assumed to be allowed to employ the \textit{existing} flexibility (offered by each user's PEV) and the diversity (resulting from the users' different usage patterns). We refer to these two as the system's elasticity. Nevertheless, the electricity consumption behaviours of the users (their PEVs' usage patterns) are not to be changed and hence the algorithm preserves users' comfort. Moreover, users' privacy is not violated as the information about their individual appliances, including PEV, is not revealed.  

The convergence criterion in Algorithm 1 can be simply assumed as a desired number of iterations of updating all users' demand profiles or it can be determined to be lower than some pre-set mean square error (MSE) between two subsequent iterations of achieving aggregated demand profiles. As we have discussed in \cite{sgtrfar}, the convergence is guaranteed to be obtained. Furthermore, users' contribution can be modelled as a cooperative game with complete information wherein a Nash equilibrium exists \cite{sgtrfar}.      

%\begin{align}
%F=\frac{\sum\limits_{n=1}^{N}(\beta_n-\alpha_n)}{\sum\limits_{n=1}^{N}\frac{E_{PEV,n}}{a_n}}.
%\end{align}

\section{Simulation Results} \label{SR}
In this section, we evaluate and present the results of our proposed model and programming technique articulated in the previous sections through computer simulations. In the simulations, we consider the number of residential users, $N$, to be 1,000 and the horizon for testing and evaluation is considered to be 24 hours for a DA programming scenario with a time granularity of one hour.

For the PEVs usage patterns, our data and distributions are based on 2009 NHTS data \cite{NHTS}. Fig. \ref{tdist} displays the distributions for arrival time, departure time, charging time and PEVs' state of charge (SOC) at the arrival time. Furthermore, we considered new standard outlets, NEMA 5-15, with 1.8 kW power transfer limit.  
We assumed that PEVs are needed to be fully charged by their respective next departure time. Additionally, we considered 24 kWh energy storage capacity for all PEVs according to Nissan Leaf model \cite{leaf}. Moreover, we adopted the PJM interconnection electricity market pricing data for both DA and RT markets in the year 2015 \cite{PJM} and assumed that PEVs are all V2G enabled.

\begin{figure}[t] 
	\centering
	\includegraphics[width=\columnwidth,height=3.6in]{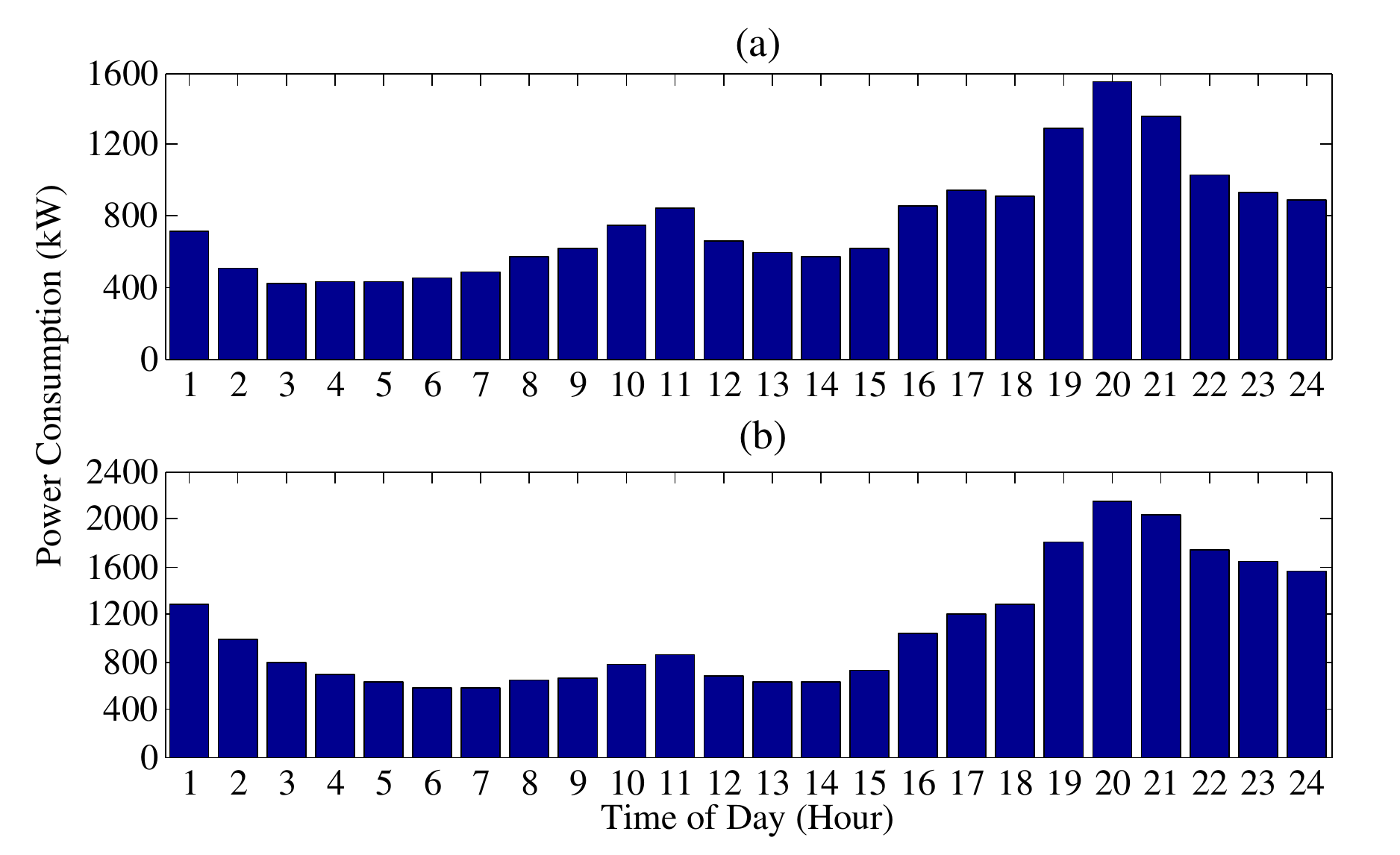}
	\caption[]{Electricity demand profile from (a) normal household appliances, i.e., without PEVs and (b) the overall electricity demand profile when users owns PEVs with different usage patterns based on NHTS data.} 
	\label{infev}\vspace{-2em}
\end{figure}

Next, we examine the DDR scheme introduced in Algorithm 1. Fig. \ref{infev} shows the assumed daily aggregated electricity demand profile of the users with and without the presence of PEVs with different usage patterns based on NHTS data.    

\begin{figure} 
	\centering
	\includegraphics[width=\columnwidth,height=2.4in]{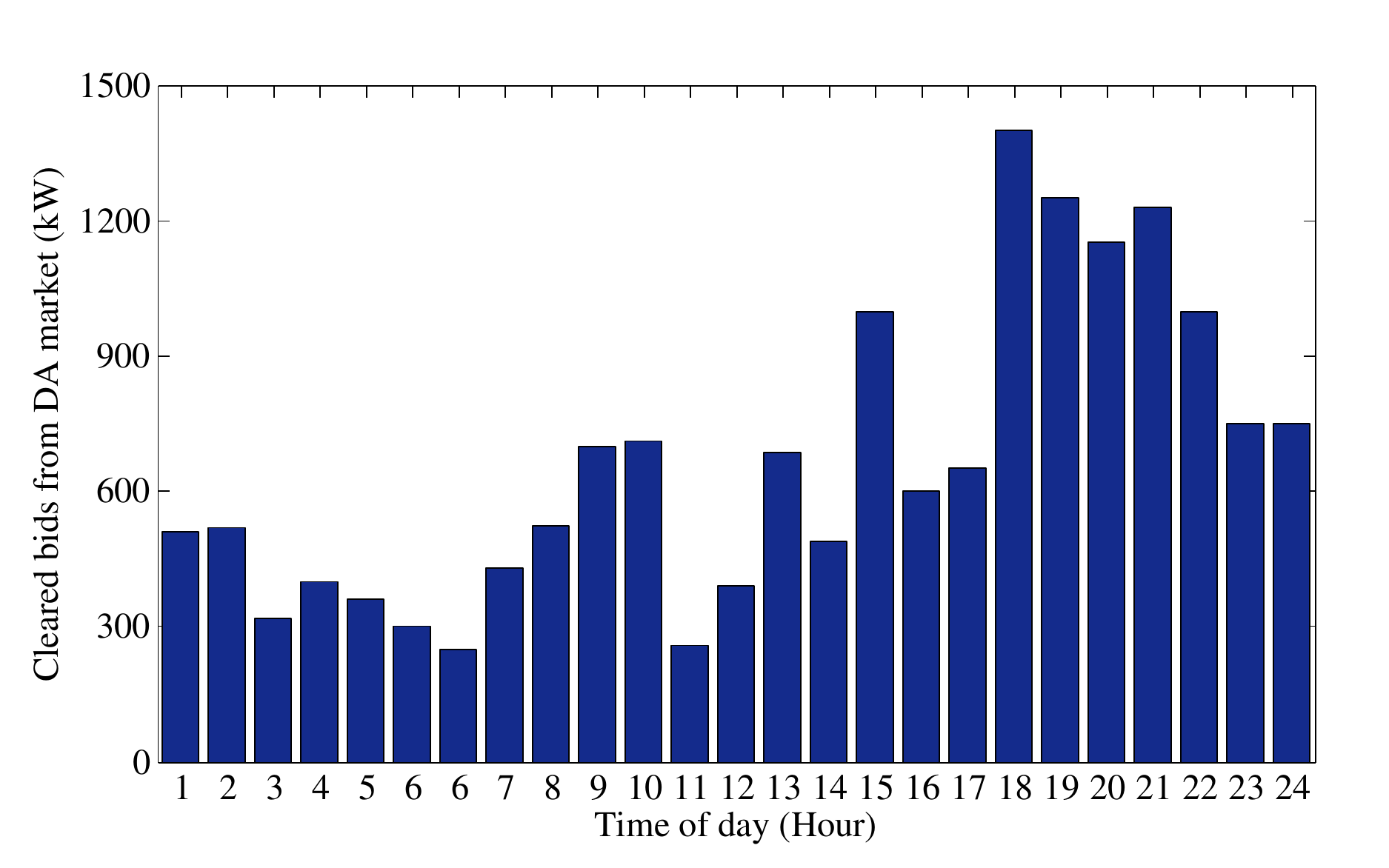}
	\caption[]{The assumed electricity profile purchased by the retailer from the DA market by the bids that could be cleared.} 
	\label{dap}\vspace{-1.7em}
\end{figure}

\begin{figure} 
	\centering
	\includegraphics[width=\columnwidth,height=2.4in]{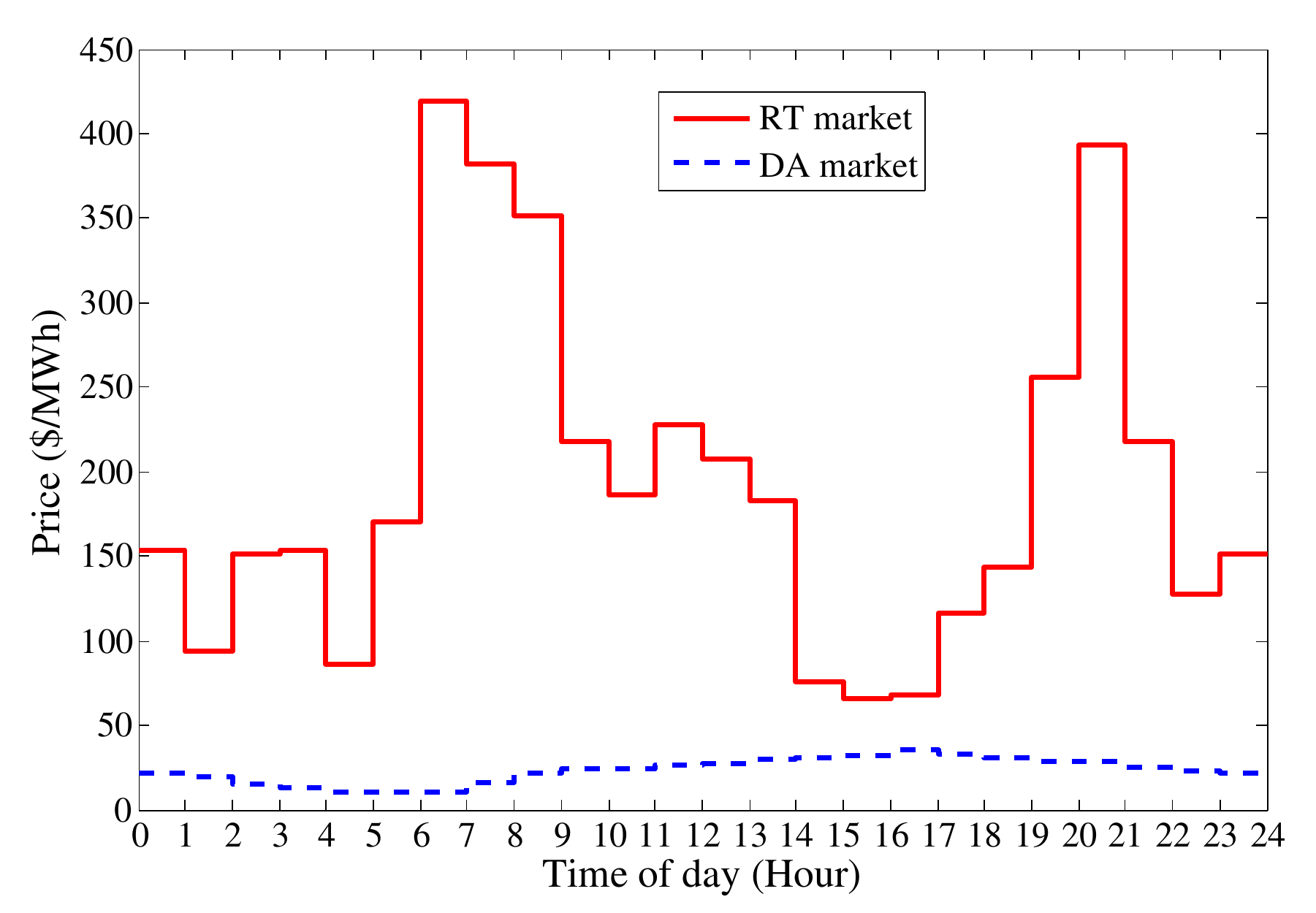}
	\caption[]{DA and RT prices for 28 November 2015 in PJM Interconnection electricity market (the number of days with such \textit{black swan} behaviour in 2014 and 2015 is several).} 
	\label{p15}\vspace{-2em}
\end{figure}

\begin{table}
	\renewcommand{\arraystretch}{2.5}
	\caption{\textlcsc{Overall energy procurement costs for the retailer}}
	\label{tab}
	\centering
	%	\resizebox{\columnwidth}{!}{  
	\begin{tabular}{|c|c|}
		\hline 
		Case  & 
		Energy procurement cost (\$)
		\\ 
		\hline
		1  &  5,382.3  \\
		\hline
		2    & 4,463.9  \\
		\hline
		3   & 4,149.2  \\
		\hline
		4   & 4,262.8 \\
		\hline
	\end{tabular}
	%	}
	\label{tab1}
\end{table}

\begin{figure}[t]
	\centering
	\includegraphics[width=\columnwidth,height=4.2in]{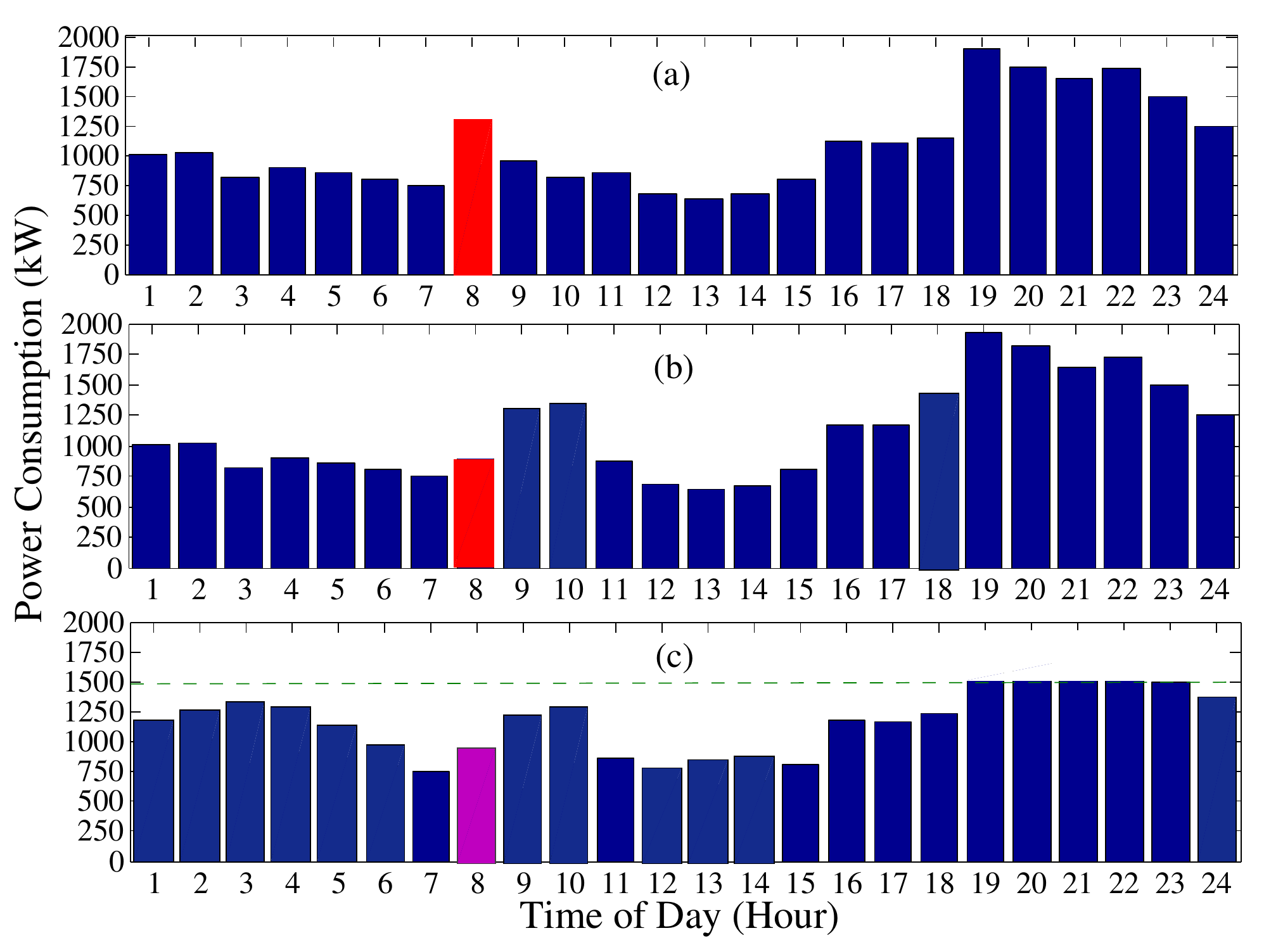}
	\caption[]{Aggregated power demand profiles for (\textit{a}) after using DR technique in \cite{sgtrfar}, (\textit{b}) when the algorithm in \cite{jointfar2} is used and (\textit{c}) when the proposed DDR algorithm is employed.}
	\label{alg}\vspace{-2em}
\end{figure}

Fig. \ref{dap} shows an assumed electricity profile cleared for the retailer in the DA market. In other words, it shows the bids that could be cleared in the market at different hours of the following day. The results of Algorithm 1 is depicted in Fig. \ref{alg}. For this, we assumed $\lambda=0.5$ in (2). For a particular day as an extreme example- Fig. \ref{p15} - it can be observed that at the eighth hour of the day the highest RT price occurs. Online demand altering can reduce the aggregated demand from 1258.7 kWh to 878.2 kWh, i.e., we can obtain almost \%30 reduction in the overall demand at that hour. This is when the constraint (5) in the proposed programming technique and algorithm is not complied. 
In case of the absence and presence of that constraint -Decarbonization constraint- DDR algorithm results can be seen in Fig. \ref{alg}. We see that the picks are curbed in Fig. 9(c). But, the power demand at eighth hour is now 985.6 kWh which cause some extra cost. This additional money could be paid back by the regulator to the retailer as subsidies for instance.   

In our simulations, convergence has been attained only after one single iteration of updating \textit{all} users' electricity demand profiles in Algorithm 1.  

It should be emphasized that this could be achieved since at that hour of the day we had almost 405 V2G enabled PEVs available at users' dwellings. Different results would be obtained for the other hours of that day. Also, it is obvious that the amount of cost savings would be dissimilar on weekdays and in the weekend.     

In Table \ref{tab1}, we compare the overall electricity procurement costs for the retailer for four cases: \textit{case (1)} is purchasing electricity without any DR, \textit{case (2)} when DR technique in \cite{sgtrfar} is used, i.e., only DA demand shaping is implemented, \textit{case (3)} when joint shaping and altering demand is applied as in \cite{jointfar2} and \textit{case (4)} when DDR is being employed. 

It can be seen that in the first case, when no DR method is used and the power demand is directly purchased from the RT market, total cost is the highest. For the second case, around \$920 is saved and in the third case the cost is further reduced by \$314.7. In the fourth case, however, when DDR is employed there is some extra cost, \$113.6, because of complying with the GHGs emissions reduction in Algorithm 1.

%\begin{figure} 
%	\centering
%	\includegraphics[width=\columnwidth,height=2.5in]{avaEV}\vspace{-0.5 em}
%	\caption[]{An PEV's expected daily power demand profile for different distributions of charging time $T$.} \vspace{-1em} 
%	\label{avaz}
%\end{figure}  

\section{Conclusion and Future Work} \label{Con}

In this paper, we proposed a fast converging and decentralized algorithm for managing V2G enabled PEVs' electricity assignments (charging and discharging) in order to simultaneously reduce the overall electricity procurement cost and GHGs emission for electricity retailers. We illustrated that with some incentives and/or regulations from the regulator, the retailer or aggregator could help lessening GHGs emissions by using our proposed decarbonized demand response (DDR) technique. 

In this work, we emphasised on the importance of considering decarbonization in DR algorithms for PEVs. However, various other combinatorial optimization methods could be investigated. Furthermore, regional GHGs emission factors can be captured into the evaluations in practice.

%\appendix
%
%\begin{proof}[Proof of proposition \ref{pro}]
%
%\end{proof}

\bibliographystyle{IEEEtran} 
\bibliography{IEEEabrv,myBIB}

\end{document}